\documentclass[12pt]{article}
\usepackage{latexsym}
\usepackage[margin=2cm,nohead]{geometry}

\newtheorem{theorem}{Theorem}[section]
\newtheorem{lemma}[theorem]{Lemma}

\def \mod {\ ({\rm mod}\,\,}
\def \deg {{\rm deg}}
\def \Nbd {{\rm Nbd}}
\def \proof {\noindent{\bf Proof}\quad}
\def \binom#1#2{{#1\choose#2}}

\newcommand{\qed}{\hfill$\Box$\vspace{0.2cm}}
\newcommand{\tfrac}[2]{{\textstyle \frac{#1}{#2}}}

\title{\bf Decomposing various graphs into short even-length cycles}

\author{
Daniel Horsley\\
School of Mathematical Sciences \\
Monash University \\
Vic 3800, Australia \\[0.1cm]
\texttt{danhorsley@gmail.com}}

\date{ }

\begin{document}
\maketitle\thispagestyle{empty}
\def\baselinestretch{1.5}\small\normalsize


\begin{abstract}
We prove that a complete bipartite graph can be decomposed into cycles of arbitrary specified
lengths provided that the obvious necessary conditions are satisfied, the length of each cycle is
at most the size of the smallest part, and the longest cycle is at most three times as long as the
second longest. We then use this result to obtain results on incomplete even cycle systems with a
hole and on decompositions of complete multipartite graphs into cycles of uniform even length.
\end{abstract}


\section{Introduction}

A decomposition of a graph $G$ is a set of subgraphs of $G$ whose edge sets partition the edge set
of $G$. In 1981 Alspach \cite{Al} posed the problem of determining when a complete graph can be
decomposed into cycles of arbitrary specified lengths. While a complete solution to this problem
continues to be elusive, a great deal of work has been done on it (see the survey \cite{Br}, for
example). Less attention has been paid to the problem of determining when a complete bipartite
graph can be decomposed into cycles of arbitrary specified lengths. Sotteau \cite{So} has solved
the problem in the case of cycles of uniform length. Chou, Fu and Huang \cite{ChFuHu} have given a
solution in the case when all the cycles have lengths in $\{4,6,8\}$. Chou and Fu \cite{ChFu} have
also examined the problem in the case when all the cycles have lengths in $\{4,m\}$, for some even
integer $m$, and have given a complete solution for $m=10$ and $m=12$. Decompositions of complete
bipartite graphs with equal sized parts into cycles and a perfect matching have also been studied,
with the case of uniform length cycles largely solved by Archdeacon et al. \cite{ArEtAl} and the
solution completed by Ma, Pu and Shen \cite{MaPuSh}. In \cite{ChFuHu} and \cite{ChFu} the analogues
of the results mentioned above for decompositions into cycles and a perfect matching are also
obtained.


Let $K_{a,b}$ denote the complete bipartite graph with parts of size $a$ and $b$, and let
$K_{a,a}-I$ denote the graph obtained from $K_{a,a}$ by removing the edges of a perfect matching.
The main results of this paper are as follows.

\begin{theorem}\label{MainTheorem}
Let $a$ and $b$ be positive integers such that either $a$ and $b$ are even or $a=b$, and let
$K^*_{a,b}$ be the graph $K_{a,b}$ if $a$ and $b$ are even and the graph $K_{a,b}-I$ if $a=b$ and
$a$ is odd. If $m_1,m_2,\ldots,m_t$ are even integers such that $4\leq m_1\leq m_2 \leq \cdots \leq
m_t \leq \min(a,b,3m_{t-1})$ and $m_1+m_2+\cdots+m_t=|E(K^*_{a,b})|$, then there is a decomposition
of $K^*_{a,b}$ into cycles of lengths $m_1,m_2,\ldots,m_t$.
\end{theorem}

Note that the conditions of Theorem \ref{MainTheorem} are necessary for the decomposition to exist,
with the exception of requiring that $m_t \leq \min(a,b,3m_{t-1})$. It is necessary that $m_t \leq
2\min(a,b)$, however.

We shall apply Theorem \ref{MainTheorem} to obtain a result concerning even cycle systems with a
hole and a result concerning decompositions of complete multipartite graphs into uniform
even-length cycles.

For positive integers $u$ and $v$ with $u \leq v$, let $K_v-K_u$ denote the graph obtained from the
complete graph of order $v$ by removing the edges of a complete subgraph of order $u$. A
decomposition of $K_v-K_u$ into $m$-cycles is often called an \emph{incomplete $m$-cycle system of
order $v$ with a hole of size $u$}. Incomplete cycle systems have received considerable attention
(see \cite{BrRo} and \cite{HoLiRo}, for example). We will concern ourselves here only with
incomplete even cycle systems. Bryant, Rodger and Spicer \cite{BrRoSp} have found necessary and
sufficient conditions for the existence of an incomplete $m$-cycle system of order $v$ with a hole
of size $u$ for all even $m$ with $m \leq 14$, and also when $v \equiv u \mod 2m)$. We prove the
following result.

\begin{theorem}\label{HoleTheorem}
Let $v$ and $u$ be odd positive integers. If $m \geq 4$ is an even integer such that
$\binom{v}{2}-\binom{u}{2} \equiv 0 \mod m)$, $u \geq m+1$, and $v-u \geq m$, then there is an
$m$-cycle decomposition of $K_v-K_u$.
\end{theorem}

Note that the condition of the theorem that $\binom{v}{2}-\binom{u}{2} \equiv 0 \mod m)$ is clearly
necessary, and also that if we demand that $u$ be the order of a non-trivial $m$-cycle system (so
we can ``fill the hole'' and obtain an $m$-cycle system of order $v$ with a subsystem of size $u$),
then the condition that $u \geq m+1$ is also necessary.

Decompositions of complete multipartite graphs into cycles have also been well studied (see, for
example, the survey \cite{Bi} and the references therein). Again, we will only consider
decompositions into cycles of uniform even length. Laskar and Auerbach \cite{LaAu} have shown that
a complete multipartite graph with all parts of equal size has a decomposition into hamilton cycles
if its vertices have even degree, and has a decomposition into hamilton cycles and a perfect
matching if its vertices have odd degree. Cavenagh and Billington \cite{CaBi} have established
necessary and sufficient conditions for the existence of a decomposition of a complete multipartite
graph into $m$-cycles for each $m \in \{4,6,8\}$. Our result is as follows.

\begin{theorem}\label{MultipartiteTheorem}
Let $G$ be a complete multipartite graph with parts of even size. If $m \geq 4$ is an even integer
such that every part of $G$ has size at least $m+2$ and $|E(G)| \equiv 0 \mod m)$, then there is a
decomposition of $G$ into $m$-cycles.
\end{theorem}

\section{Preliminary Results} \label{PreliminaryResultsSection}

Our main goal in Sections \ref{PreliminaryResultsSection}, \ref{JoinCyclesSection} and
\ref{BipartiteDecompSection} is to prove Theorem \ref{MainTheorem}. In this section our aim is to
prove Lemmas \ref{PathSwitchGeneral}, \ref{PathThing} and \ref{FigureOfEight}, which we will
require in Section \ref{JoinCyclesSection}. Lemma \ref{FigureOfEightHelper} is used only in the
proof of Lemma \ref{FigureOfEight}.

We require some more notation before we proceed. The \emph{size} of a graph is the number of edges
it contains. We say that a graph is \emph{even} if each of its vertices has even degree and that it
is \emph{odd} if each of its vertices has odd degree. The \emph{neighbourhood} of a vertex $x$ in a
graph $G$, denoted $\Nbd_G(x)$, is the set of all the vertices adjacent in $G$ to $x$ (not
including $x$ itself). If $M = m_1,m_2,\ldots,m_t$ is a list of integers and $G$ is a graph which
is either even or odd, then an \emph{$(M)^*$-packing} of $G$ is a decomposition of some subgraph of
$G$ into $t$ cycles of lengths $m_1,m_2,\ldots,m_t$ if $G$ is an even graph, and into $t$ cycles of
lengths $m_1,m_2,\ldots,m_t$ and a perfect matching if $G$ is an odd graph. The \emph{leave} of an
$(M)^*$-packing of $G$ is the complement in $G$ of the subgraph of $G$ which is decomposed. In
other words, the leave is the spanning subgraph of $G$ whose edges are precisely those which do not
appear in the cycles or the perfect matching (if it exists) of the packing. Note that the leave of
an $(M)^*$-packing is always an even graph. An $(M)^*$-decomposition of $G$ is an $(M)^*$-packing
of $G$ whose leave is an empty graph. If we know $G$ to be an even graph, then we will drop the
asterisk and simply refer to an $(M)$-packing or $(M)$-decomposition of $G$.

The {\em length} of a cycle or path is the number of edges it contains. A cycle of length $p$ is
called a $p$-cycle and a path of length $q$ is called a $q$-path. The $p$-cycle with vertices
$x_1,x_2,\ldots,x_p$ and edges $x_1x_2,x_2x_3,\ldots,x_{p-1}x_p,x_px_1$ is denoted by
$(x_1,x_2,\ldots,x_p)$ and the $q$-path with vertices $y_0,y_1,\ldots,y_q$ and edges
$y_0y_1,y_1y_2,\ldots,y_{q-1}y_q$ is denoted by $[y_0,y_1,\ldots,y_q]$. We will say that $y_0$ to
$y_q$ are the \emph{end vertices} of such a path, and that the path is from $y_0$ to $y_q$ (or from
$y_q$ to $y_0$). For graphs $G$ and $H$, the union of $G$ and $H$, denoted by $G \cup H$, is the
graph with vertex set $V(G) \cup V(H)$ and edge set $E(G) \cup E(H)$.

We define two special kinds of graphs which we call rings and chains. Note that these definitions
vary slightly from the definitions given in \cite{BrHoLong} to suit our purposes in this paper.

A graph $G$ is an {\em $(a_1,a_2,\ldots,a_s)$-chain} if $G$ is the edge-disjoint union of $s\geq 2$
cycles $A_1, A_2, \ldots, A_s$ such that
\begin{itemize}
\item $A_i$ is a cycle of length $a_i$ for $i=1,2,\ldots,s$; and
\item for $1\leq i<j\leq s$, $|V(A_i)\cap V(A_j)|=1$ if $j=i+1$ and
$|V(A_i)\cap V(A_j)|=0$ otherwise.
\end{itemize}
We call $A_1$ and $A_s$ the \emph{end cycles} of $G$. A graph is an {\em $s$-chain}, or just a {\em
chain}, if it is an $(a_1,a_2,\ldots,a_s)$-chain for some integers $a_1,a_2,\ldots,a_s$. We denote
a $(p,q)$-chain with cycles $(x_1,x_2,\ldots,x_{p-1},c)$ and $(y_1,y_2,\ldots,y_{q-1},c)$ by
$(x_1,x_2,\ldots,x_{p-1},c)\cdot(c,y_1,y_2,\ldots,y_{q-1})$.

A graph $G$ is an {\em $(a_1,a_2,\ldots,a_s)$-ring} if $G$ is the edge-disjoint union of $s\geq 2$
cycles $A_1, A_2, \ldots, A_s$ such that
\begin{itemize}
\item $A_i$ is a cycle of length $a_i$ for $i=1,2,\ldots,s$;
\item if $s \geq 3$, then for $1\leq i<j\leq s$, $|V(A_i)\cap V(A_j)|=1$ if $j=i+1$ or if $(i,j)=(1,s)$, and
$|V(A_i)\cap V(A_j)|=0$ otherwise; and
\item if $s=2$, then $|V(A_1)\cap V(A_2)|=2$.
\end{itemize}
We call $A_1, A_2, \ldots, A_s$ the \emph{ring cycles} of $G$ (note that $G$ contains cycles which
are not ring cycles). A graph is an {\em $s$-ring}, or just a {\em ring}, if it is an
$(a_1,a_2,\ldots,a_s)$-ring for some integers $a_1,a_2,\ldots,a_s$.

The techniques used in this paper are often the bipartite analogues of techniques used for complete
graphs in \cite{BrHoLong} and \cite{BrHoBigN}. In particular, the crucial Lemma \ref{JoinCycles} of
this paper is the bipartite analogue of Lemma 3.1 of \cite{BrHoBigN}. Unfortunately however, we
have thus far been unable to find a good bipartite analogue for Lemma 1.2 of \cite{BrHoLong}. We
begin by using Lemma 2.1 of \cite{BrHoLong}, which deals with modifying packings of complete
graphs, to obtain a similar result for decompositions of general graphs. Given a permutation $\pi$
of a set $V$ and a graph $G$ with $V(G)\subseteq V$, the graph $\pi(G)$ is defined to be the graph
with vertex set $\{\pi(x):x\in V(G)\}$ and edge set $\{\pi(x)\pi(y):xy\in E(G)\}$.

\begin{lemma}\label{PathSwitchGeneral}
Let $G$ be a graph, let $M$ be a list of integers, let $\mathcal{P}$ be an $(M)^*$-packing of $G$,
let $L$ be the leave of $\mathcal{P}$, and let $\alpha$ and $\beta$ be vertices in $G$ such that
$\Nbd_G(\alpha)=\Nbd_G(\beta)$, and let $\pi$ be the transposition $(\alpha \beta)$. Then there
exists a partition of the set $(\Nbd_L(\alpha)\cup \Nbd_L(\beta)) \setminus ((\Nbd_L(\alpha)\cap
\Nbd_L(\beta)) \cup \{\alpha,\beta\})$ into pairs such that for each pair $\{u,v\}$ of the
partition, there exists an $(M)^*$-packing $\mathcal{P}'$ of $G$ whose leave $L'$ differs from $L$
only in that $\alpha u$, $\alpha v$, $\beta u$ and $\beta v$ are edges in $L'$ if and only if they
are not edges in $L$. Furthermore, if $\mathcal{P} = \{C_1,C_2,\ldots,C_t\}$ ($G$ even) or
$\mathcal{P} = \{F,C_1,C_2,\ldots,C_t\}$ ($G$ odd) where $C_1,C_2,\ldots,C_t$ are cycles and $F$ is
a perfect matching, then $\mathcal{P}' = \{C'_1,C'_2,\ldots,C'_t\}$ ($G$ even) or $\mathcal{P}' =
\{F',C'_1,C'_2,\ldots,C'_t\}$ ($G$ odd), where, for $i=1,2,\ldots,t$, $C'_i$ is a cycle of the same
length as $C_i$ and $F'$ is a perfect matching such that
\begin{itemize}
    \item[(i)]
either $F'=F$ or $F'=\pi(F)$;
    \item[(ii)]
for $i = 1,2,\ldots,t$ if neither $\alpha$ nor $\beta$ is in $V(C_i)$, then $C'_i=C_i$;
    \item[(iii)]
for $i = 1,2,\ldots,t$ if exactly one of $\alpha$ and $\beta$ is in $V(C_i)$, then either
$C'_i=C_i$ or $C'_i=\pi(C_i)$; and
    \item[(iv)]
for $i = 1,2,\ldots,t$ if both $\alpha$ and $\beta$ are in $V(C_i)$, then $C'_i \in
\{C_i,\pi(C_i),\pi(P_i)\cup P^{\dag}_i,P_i \cup \pi(P^{\dag}_i)\}$, where $P_i$ and
$P^{\dag}_i$ are the two paths in $C_i$ which have end vertices $\alpha$ and $\beta$.
\end{itemize}
\end{lemma}

\proof If $G$ is an even graph and $|V(G)|$ is even let $G^{\dag}$ be the graph obtained from $G$
by adding an isolated vertex. Otherwise, let $G^{\dag}=G$. Note that we can consider $\mathcal{P}$
as an $(M)^*$-packing of the complete graph on vertex set $V(G^{\dag})$. The leave $H$ of this
packing is the edge-disjoint union of $L$ and the complement of $G^{\dag}$.  By applying Lemma 2.1
of \cite{BrHoLong} we can obtain a partition of the set $(\Nbd_H(\alpha)\cup \Nbd_H(\beta))
\setminus ((\Nbd_H(\alpha)\cap \Nbd_H(\beta)) \cup \{\alpha,\beta\})$ into pairs such that for each
pair $\{u,v\}$ of the partition, there exists an $(M)^*$-packing $\mathcal{P}'$ of the complete
graph on the vertex set  $V(G^{\dag})$ whose leave $H'$ differs from $H$ only in that $\alpha u$,
$\alpha v$, $\beta u$ and $\beta v$ are edges in $H$ if and only if they are not edges in $H$.
Furthermore, if $\mathcal{P} = \{C_1,C_2,\ldots,C_t\}$ ($|V(G^{\dag})|$ even) or $\mathcal{P} =
\{F,C_1,C_2,\ldots,C_t\}$ ($|V(G^{\dag})|$ odd) where $C_1,C_2,\ldots,C_t$ are cycles and $F$ is a
perfect matching, then $\mathcal{P}' = \{C'_1,C'_2,\ldots,C'_t\}$ ($|V(G^{\dag})|$ even) or
$\mathcal{P}' = \{F',C'_1,C'_2,\ldots,C'_t\}$ ($|V(G^{\dag})|$ odd) where for $i=1,2,\ldots,t$,
$C'_i$ is a cycle of the same length as $C_i$ and $F'$ is a perfect matching such that (i), (ii),
(iii) and (iv) hold. Since $\Nbd_G(\alpha)=\Nbd_G(\beta)$ and $H$ is the edge-disjoint union of $L$
and the complement of $G^{\dag}$, we have that $(\Nbd_H(\alpha)\cup \Nbd_H(\beta)) \setminus
((\Nbd_H(\alpha)\cap \Nbd_H(\beta) \cup \{\alpha,\beta\}) = (\Nbd_L(\alpha)\cup \Nbd_L(\beta))
\setminus ((\Nbd_L(\alpha)\cap \Nbd_L(\beta)) \cup \{\alpha,\beta\})$. Furthermore, $\alpha u,
\alpha v, \beta u, \beta v \in E(G)$, so $H'$ is the edge-disjoint union of the complement of
$G^{\dag}$ and a graph $L'$ on vertex set $V(G)$ which differs from $L$ only in that $\alpha u$,
$\alpha v$, $\beta u$ and $\beta v$ are edges in $L'$ if and only if they are not edges in $L$.
Thus, by viewing $\mathcal{P}'$ as a packing of $G$, we have our result. \qed

\noindent{\bf Remark.} In this paper we will always use the above result in the case where $G$ is a
complete bipartite graph. Thus the condition $\Nbd_G(\alpha)=\Nbd_G(\beta)$ will hold if and only
if $\alpha$ and $\beta$ are in the same part of $G$.

When applying Lemma \ref{PathSwitchGeneral}, we will say that $\mathcal{P}'$ is the {\em
$(M)^*$-packing of $K_{a,b}$ obtained from $\mathcal{P}$ by performing the $(\alpha,\beta)$-switch
with origin $u$ and terminus $v$} (we could equally call $v$ the origin and $u$ the terminus).

In the remainder of this paper we will often prove results relating to an $(M)^*$-packing of
$K_{a,b}$, where we take it as read that $M$ is a list of positive even integers and that $a$ and
$b$ are positive integers such that either both $a$ and $b$ are even or $a=b$. Note that when $a$
and $b$ are even this is a packing of $K_{a,b}$ with cycles of lengths $M$ and when $a=b$ is odd
this is a packing of $K_{a,a}$ with cycles of lengths $M$ and a perfect matching.

\begin{lemma}\label{PathThing}
Suppose that there exists an $(M)^*$-packing $\mathcal{P}$ of $K_{a,b}$ with a leave of size $l$
whose only non-trivial component $H$ contains a path $P=[x_0,x_1,\ldots,x_t]$ of even length at
least $4$ such that the edges in $E(H) \setminus E(P)$ form a path and such that $x_1x_t \notin
E(H)$. Let $S$ be the $(x_0,x_t)$-switch with origin $x_1$ (note that $x_0$ and $x_t$ are in the
same part of $K_{a,b}$ since $P$ is a path of even length) and let $\mathcal{P}'$ be the
$(M)^*$-packing of $K_{a,b}$ obtained from $\mathcal{P}$ by performing $S$. If $S$ does not have
terminus $x_{t-1}$, then the leave of $\mathcal{P}'$ has a decomposition into a $t$-cycle and an
$(l-t)$-cycle, and there are at least as many vertices of degree $4$ in the leave of $\mathcal{P}'$
as there are in the leave of $\mathcal{P}$.
\end{lemma}

\noindent {\bf Remark.} In particular, notice that if $H$ contains a vertex of degree at least $4$,
then the leave of $\mathcal{P}'$ will have only one non-trivial component.

\proof Since each vertex of $H$ has even degree it is clear that $\deg_H(x_0)=\deg_H(x_t)=2$. Let
the path with edge set $E(H) \setminus E(P)$ be $[x_0=y_0,y_1,y_2,\ldots,y_{l-t}=x_t]$. Since $S$
does not have terminus $x_{t-1}$, it must have terminus $y_1$ or $y_{l-t-1}$. If $S$ has terminus
$y_1$, then the leave of $\mathcal{P}'$ has a decomposition into the $t$-cycle
$(x_1,x_2,\ldots,x_t)$ and the $(l-t)$-cycle $(y_1,y_2,\ldots,y_{l-t})$, and there is one more
vertex of degree $4$ in the leave of $\mathcal{P}'$ than there is in the leave of $\mathcal{P}$.
Otherwise $S$ has terminus $y_{l-t-1}$, the leave of $\mathcal{P}'$ has a decomposition into the
$t$-cycle $(x_1,x_2,\ldots,x_t)$ and the $(l-t)$-cycle $(y_0,y_1,y_2,\ldots,y_{l-t-1})$, and there
are the same number of vertices of degree $4$ in the leave of $\mathcal{P}'$ as there are in the
leave of $\mathcal{P}$. \qed

\begin{lemma}\label{FigureOfEightHelper}
Suppose there exists an $(M)^*$-packing of $K_{a,b}$ with a leave whose only non-trivial component
is a $(p,q)$-chain. If $m$ is an even integer such that $p\leq m$ and $p+q-m\geq 4$, then there
exist
\begin{itemize}
    \item[(i)]
an $(M)^*$-packing of $K_{a,b}$ with a leave whose only non-trivial component either has a
decomposition into an $m$-cycle and a $(p+q-m)$-cycle, or is an $(m-p+2,2p+q-m-2)$-chain; and
    \item[(ii)]
an $(M)^*$-packing of $K_{a,b}$ with a leave whose only non-trivial component either has a
decomposition into an $m$-cycle and a $(p+q-m)$-cycle, or is an $(m-p+4,2p+q-m-4)$-chain.
\end{itemize}
\end{lemma}

\proof If $p=m$ then we are finished, so assume $p \leq m-2$. Let $\mathcal P$ be an
$(M)^*$-packing of $K_{a,b}$ and let $H=(x_1,x_2,\ldots,x_{p-1},c)\cdot(c,y_1,y_2,\ldots,y_{q-1})$
be the only non-trivial component in its leave.

{\bf (i)} Let $\mathcal{P}'$ be the $(M)^*$-packing of $K_{a,b}$ obtained from $\mathcal P$ by
performing the $(x_1,y_{m-p+1})$-switch $S'$ with origin $x_2$ (note that $x_1$ and $y_{m-p+1}$ are
in the same part of $K_{a,b}$ since $[x_1,c,y_1,y_2,\ldots,y_{m-p+1}]$ is an $(m-p+2)$-path and
$m-p+2$ is even). If $S'$ has terminus $y_{m-p}$, then the only non-trivial component in the leave
of $\mathcal P'$ is
$$(y_1,y_2,\ldots,y_{m-p},x_1,c)\cdot
(c,x_{p-1},x_{p-2},\ldots,x_2,y_{m-p+1},y_{m-p+2},\ldots,y_{q-1}),$$ which is an
$(m-p+2,2p+q-m-2)$-chain. Otherwise $S'$ does not have terminus $y_{m-p}$ and, by Lemma
\ref{PathThing} with $P=[x_1,x_2,\ldots,x_{p-1},c,y_1,y_2,\ldots,y_{m-p+1}]$, the leave of
$\mathcal P'$ has a decomposition into an $m$-cycle and a $(p+q-m)$-cycle.

{\bf (ii)} Let $\mathcal{P}''$ be the $(M)^*$-packing of $K_{a,b}$ obtained from $\mathcal{P}$ by
performing the $(x_2,y_{m-p+2})$-switch $S''$ with origin $x_3$ (note that $x_2$ and $y_{m-p+2}$
are in the same part of $K_{a,b}$ since $[x_2,x_1,c,y_1,y_2,\ldots,y_{m-p+2}]$ is an $(m-p+4)$-path
and $m-p+4$ is even). If $S''$ has terminus $y_{m-p+1}$, then the only non-trivial component in the
leave of $\mathcal{P}''$ is
$$(y_1,y_2,\ldots,y_{m-p+1},x_2,x_1,c)\cdot
(c,x_{p-1},x_{p-2},\ldots,x_3,y_{m-p+2},y_{m-p+3},\ldots,y_{q-1}),$$ which is an
$(m-p+4,2p+q-m-4)$-chain. Otherwise, $S''$ does not have terminus $y_{m-p+1}$ and, by Lemma
\ref{PathThing} with $P=[x_2,x_3,\ldots,x_{p-1},c,y_1,y_2,\ldots,y_{m-p+2}]$, the leave of
$\mathcal{P}''$ has a decomposition into an $m$-cycle and a $(p+q-m)$-cycle. \qed

\begin{lemma}\label{FigureOfEight}
Suppose there exists an $(M)^*$-packing of $K_{a,b}$ with a leave whose only non-trivial component
is a $(p,q)$-chain. If $m_1$ and $m_2$ are even integers such that $m_1,m_2\geq 4$ and
$m_1+m_2=p+q$, then there exists an $(M,m_1,m_2)^*$-decomposition of $K_{a,b}$.
\end{lemma}

\proof We can assume without loss of generality that $p \leq q$ and $m_1 \leq m_2$. Since
$m_1+m_2=p+q$, this implies that $p \leq m_2$. Let $m=m_2$. We wish to find an
$(M,m,p+q-m)^*$-decomposition of $K_{a,b}$.

If $p=m$ then we are finished immediately, so assume $m>p$. We apply Lemma
\ref{FigureOfEightHelper} (i) if $p\geq\frac{m+4}2$ and we apply Lemma \ref{FigureOfEightHelper}
(ii) if $p\leq\frac{m+2}2$. Thus, we either obtain the required $(M)^*$-packing of $K_{a,b}$ or we
obtain an $(M)^*$-packing of $K_{a,b}$ with a leave whose only non-trivial component is an
$(m-p+2,2p+q-m-2)$-chain when $p\geq\frac{m+4}2$, and is an $(m-p+4,2p+q-m-4)$-chain when
$p\leq\frac{m+2}2$. Let $p_2=m-p+2$ if $p\geq\frac{m+4}2$ and let $p_2=m-p+4$ if
$p\leq\frac{m+2}2$. If $p_2=m$ we are finished. Otherwise, we now apply Lemma
\ref{FigureOfEightHelper} (i) if $p_2\geq\frac{m+4}2$ and we apply Lemma \ref{FigureOfEightHelper}
(ii) if $p_2\leq\frac{m+2}2$. We claim that by repeating this process, in each instance applying
Lemma \ref{FigureOfEightHelper} (i) when $p_i\geq\frac{m+4}2$, applying Lemma
\ref{FigureOfEightHelper} (ii) when $p_i\leq\frac{m+2}2$, and defining $p_{i+1}$ by
$p_{i+1}=m-p_i+2$ if $p_i\geq\frac{m+4}2$ and $p_{i+1}=m-p_i+4$ if $p_i\leq\frac{m+2}2$, we
eventually obtain the required $(M)^*$-packing of $K_{a,b}$. To see this, observe that if $m\equiv
0 \mod 4)$ and $p=\frac{m+4}2$, then the sequence $p,p_2,p_3,\ldots$ is
$$\textstyle\frac{m+4}2,\frac m2,\frac{m+8}2,\frac{m-4}2,\ldots,m-4,6,m-2,4,m$$
and that if $m\equiv 2\mod 4)$ and $p=\frac{m+2}2$, then the sequence $p,p_2,p_3,\ldots$ is
$$\textstyle\frac{m+2}2,\frac{m+6}2,\frac{m-2}2,\frac{m+10}2,\ldots,m-4,6,m-2,4,m.$$
In either case, the sequence contains every even integer $x$ in the range $4\leq x\leq m$. \qed

\section{Merging Cycles}\label{JoinCyclesSection}

Our main goal in this section is to prove Lemma \ref{JoinCycles} which allows us to merge a number
of cycles in a decomposition into one large cycle, given the existence of a ``catalyst cycle'' of
appropriate length. We require a number of preliminary lemmas first, however, many of which will
also be used in Sections 5 and 6.

\begin{lemma}\label{PathOfCycles}
Suppose there exists an $(M)^*$-packing of $K_{a,b}$, where $a \leq b$, with a leave of size $l$,
where $l \leq 2a+2$ if $a < b$ and $l \leq 2a$ if $a = b$, with only one non-trivial component $H$.
If $m_1$ and $m_2$ are even integers such that $m_1,m_2 \geq 4$ and either
\begin{itemize}
    \item
$H$ is a chain which has a decomposition into an $m_1$-path and an $m_2$-path; or
    \item
$H$ is a ring which has a decomposition into an $m_1$-path and an $m_2$-path;
\end{itemize}
then there exists an $(M,m_1,m_2)^*$-decomposition of $K_{a,b}$.
\end{lemma}

\proof Let $\mathcal{P}$ be an $(M)^*$-packing of $K_{a,b}$ which satisfies the hypotheses of the
lemma and let $L$ be its leave. We proceed by induction. First we show that the lemma is true when
$H$ is a $2$-chain or a $2$-ring.

If $H$ is a $2$-chain then $\mathcal{P}$ satisfies the hypotheses of Lemma \ref{FigureOfEight} and
we are finished, so we may suppose that $H$ is a $2$-ring. Since $l \leq 2a+2$ if $a < b$ and $l
\leq 2a$ if $a = b$, there must exist vertices $u$ and $v$ in the same part of $K_{a,b}$ such that
$\deg_L(u) = 4$ and $\deg_L(v) = 0$ (to see this, separately consider two cases depending on
whether the two vertices of degree $4$ in $L$ are in the same part of $K_{a,b}$). Let $y$ be a
neighbour in $L$ of $u$. Let $\mathcal{P}^{\dag}$ be the $(M)^*$-packing of $K_{a,b}$ obtained from
$\mathcal{P}$ by performing the $(u,v)$-switch $S$ with origin $y$. Regardless of the terminus of
$S$, the only non-trivial component of the leave of $\mathcal{P}$ is a $2$-chain and we can apply
Lemma \ref{FigureOfEight} to complete the proof.

We now show that, for each integer $s \geq 3$, if the lemma holds when $H$ is an $(s-1)$-chain or
an $(s-1)$-ring, then it holds when $H$ is an $s$-chain or an $s$-ring. The proof splits into two
cases.

{\bf Case 1.} Suppose that $H$ is an $s$-chain. Let $P=[x_0,x_1,\ldots,x_{m_1}]$ be an $m_1$-path
in $H$ such that the edges in $E(H) \setminus E(P)$ form an $m_2$-path. Since every vertex of
degree $4$ in $H$ must be a vertex of $P$ which is not an end vertex, it is easy to see that $x_0$
and $x_{m_1}$ are vertices of degree $2$ in $H$ such that $x_0$ is in one end cycle of $H$ and
$x_{m_1}$ is in the other end cycle of $H$. Note that $x_0$ and $x_{m_1}$ are in the same part of
$K_{a,b}$ since $m_1$ is even. Let $L^{\dag}$ be the leave of the $(M)^*$-packing of $K_{a,b}$
obtained from $\mathcal{P}$ by performing the $(x_0,x_{m_1})$-switch $S$ with origin $x_1$. If $S$
does not have terminus $x_{m_1-1}$, then, by Lemma \ref{PathThing}, $L^{\dag}$ has a decomposition
into an $m_1$-cycle and an $m_2$-cycle, and we are finished (note that $x_1x_{m_1} \notin E(H)$
since $H$ is an $s$-chain with $s \geq 3$). Otherwise $S$ has terminus $x_{m_1-1}$, and it can be
seen that the only non-trivial component $H^{\dag}$ of $L^{\dag}$ is an $(s-1)$-ring. Furthermore,
if we let $P^{\dag}$ be the path with $E(P^{\dag})= (E(P) \setminus \{x_0x_1,x_{m_1-1}x_{m_1}\})
\cup \{x_0x_{m_1-1},x_1x_{m_1}\}$, then it is clear that $P^{\dag}$ is an $m_1$-path in $H^{\dag}$
such that the edges in $E(H^{\dag}) \setminus E(P^{\dag})$ form an $m_2$-path. Thus we are finished
by the inductive hypothesis.

{\bf Case 2.} Suppose that $H$ is an $s$-ring. Let $P=[x_0,x_1,\ldots,x_{m_1}]$ be an $m_1$-path in
$H$ such that the edges in $E(H) \setminus E(P)$ form an $m_2$-path. Since every vertex of degree
$4$ in $H$ must be a vertex of $P$ which is not an end vertex, it is easy to see that $x_0$ and
$x_{m_1}$ are vertices of degree $2$ in $H$ and that there is a ring cycle $C$ of $H$ such that
$x_0, x_{m_1} \in V(C)$. Let $u$ be a vertex of $C$ such that $\deg_H(u)=4$ and such that there is
at least one other vertex of degree $4$ in $L$ in the same part of $K_{a,b}$ as $u$ (such a vertex
$u$ must exist since $s \geq 3$ and hence there are at least three vertices of degree $4$ in $L$).
Let $v$ be a vertex in the same part of $K_{a,b}$ as $u$ such that $\deg_L(v)=0$ (such a vertex
must exist since $l \leq 2a+2$). Let $y$ and $z$ be the neighbours in $C$ of $u$. Let $D$ be the
ring cycle of $H$, other than $C$, such that $u \in V(D)$. Let $L^{\dag}$ be the leave of the
$(M)^*$-packing of $K_{a,b}$ obtained from $\mathcal{P}$ by performing the $(u,v)$-switch $S$ with
origin $y$. If the terminus of $S$ is a neighbour in $D$ of $u$, then the only non-trivial
component $H^{\dag}$ of $L^{\dag}$ is an $(s-1)$-ring and it can be seen that $H^{\dag}$ has a
decomposition into an $m_1$-path and an $m_2$-path (take an $m_1$-path with edge set $(E(P)
\setminus (E(C) \cup E(D))) \cup E(P') \cup E(P'')$ where $P'$ and $P''$ are paths in the ring
cycle of $H^{\dag}$ which is not a ring cycle of $H$ whose lengths add to $|E(P) \cap (E(C) \cup
E(D))|$). Thus we are finished by the inductive hypothesis. Otherwise, the terminus of $S$ is $z$,
the only non-trivial component $H^{\dag}$ of $L^{\dag}$ is an $s$-chain, and it can be seen that
$H^{\dag}$ has a decomposition into an $m_1$-path and an $m_2$-path (take an $m_1$-path with edge
set $(E(P) \setminus E(C)) \cup E(P')$ where $P'$ is a path of length $|E(P) \cap E(C)|$ in the
cycle of $H^{\dag}$ which is not a ring cycle of $H$). Thus we can proceed as we did in Case 1.
\qed

\begin{lemma}\label{AttachCycles}
Suppose that there exists an $(M)^*$-packing of $K_{a,b}$ with a leave $L$ of size $l$ with $k$
non-trivial components such that exactly one vertex of $L$ has degree $4$ and every other vertex of
$L$ has degree $2$ or degree $0$. If $R$ is one of the parts of $K_{a,b}$ and $m_1$ and $m_2$ are
integers such that $m_1,m_2 \geq k+1$ and $m_1+m_2=l$, then there exists $(M)^*$-packing of
$K_{a,b}$ with a leave whose only non-trivial component is a chain which has a decomposition into
an $m_1$-path and an $m_2$-path such that if $m_1,m_2 \geq 3$ then at least one end vertex of the
paths is in $R$.
\end{lemma}

\proof Let $\mathcal{P}$ be an $(M)^*$-packing of $K_{a,b}$ which satisfies the hypotheses of the
lemma. It follows that $L$ has one component which is a $2$-chain and that any other non-trivial
component of $L$ is a cycle. Let the non-trivial components of $L$ be $H$ and
$C_1,C_2\ldots,C_{k-1}$ (just $H$ if $k=1$), where $H$ is the $2$-chain. Let $h$ be the size of $H$
and, if $k \geq 2$, let $c_i$ be the size of $C_i$ for $i \in \{1,2,\ldots,k-1\}$.

For each integer $s$ such that $1 \leq s \leq k$ we claim that, for all integers $x_1$ and $x_2$
such that $x_1 \geq s+1$, $x_2 \geq s+1$ and $x_1+x_2=h+c_1+c_2+\cdots+c_{s-1}$ ($x_1+x_2=h$ if
$s=1$), there exists an $(M)^*$-packing of $K_{a,b}$ with a leave $L''$ such that the non-trivial
components of $L''$ are $H''$ and $C_{s},C_{s+1},\ldots,C_{k-1}$ (just $H''$ if $s=k$), where $H''$
is a chain which has a decomposition into an $x_1$-path and an $x_2$-path such that if $x_1,x_2
\geq 3$ then at least one end vertex of the paths is in $R$. Note that if this claim holds for
$s=k$ then we have the required result, so it suffices to prove the claim. We will do this by
induction on $s$.

Suppose first that $s=1$. For any integers $x_1$ and $x_2$ such that $x_1,x_2 \geq 2$ and
$x_1+x_2=h$, it is easy to see that $H$ has a decomposition into an $x_1$-path and an $x_2$-path
such that if $x_1,x_2 \geq 3$ then at least one end vertex of the paths is in $R$. Thus our claim
is true when $s=1$.

Now suppose that our claim holds when $s=t$ for some integer $t$ such that $1 \leq t < k$. We will
show that it holds when $s=t+1$. Let $x_1$ and $x_2$ be integers such that $x_1,x_2 \geq t+2$ and
$x_1+x_2=h+c_1+c_2+\cdots+c_t$. Suppose without loss of generality that $x_1 \leq x_2$. Then we
have that $t+2 \leq x_1 \leq \frac{1}{2}(h+c_1+c_2+\cdots+c_t)$. It is easy to see that we can find
positive integers $p$ and $p^{\dag}$ such that $x_1=p+p^{\dag}$, $t+1 \leq p \leq
\frac{1}{2}(h+c_1+c_2+\cdots+c_{t-1})$ and $1 \leq p^{\dag} \leq \frac{1}{2}c_t$. Furthermore, we
can choose these integers such that $p \geq 3$ if $x_1 \geq 4$. Let $q=h+c_1+c_2+\cdots+c_{t-1}-p$
and note that $q \geq p$ since $p \leq \frac{1}{2}(h+c_1+c_2+\cdots+c_{t-1})$.

Thus, by our inductive hypothesis, there exists an $(M)^*$-packing $\mathcal{P}'$ of $K_{a,b}$ with
a leave $L'$ such that the non-trivial components of $L'$ are $H'$ and
$C_{t},C_{t+1},\ldots,C_{k-1}$, where $H'$ is an $f$-chain, for some integer $f \geq 2$, which has
a decomposition into a $p$-path $P$ and a $q$-path $Q$ such that if $p \geq 3$ then at least one
end vertex of $P$ and $Q$ is in $R$. Let $u$ and $w$ be the end vertices of $P$ and $Q$, where $w
\in R$ if $p \geq 3$.

Let $y$ be the neighbour in $P$ of $u$, let $z$ be the neighbour in $Q$ of $u$, and let $v$ be a
vertex in $C_t$ which is in the same part of $K_{a,b}$ as $u$. Let $L''$ be the leave of the
$(M)^*$-packing of $K_{a,b}$ obtained from $\mathcal{P}'$ by performing the $(u,v)$-switch $S$ with
origin $y$. Then the non-trivial components of $L''$ are $H''$ and $C_{t+1},C_{t+2},\ldots,C_{k-1}$
(or just $H''$ if $t=k-1$) where $H''$ is an $(f+1)$-chain if the terminus of $S$ is $z$ and $H''$
is an $f$-chain otherwise (if the terminus of $S$ is not $z$ then it is a neighbour in $C_t$ of
$v$). In either case it can be seen that $H''$ has a decomposition into an $x_1$-path and an
$x_2$-path such that $w$ is an end vertex of each path (take an $x_1$-path with edge set $(E(P)
\setminus \{uy\}) \cup \{vy\} \cup E(P^{\dag})$ where $P^{\dag}$ is a path of length $p^{\dag}$
whose edges are in $E(C_t) \cap E(H'')$). The proof is completed by noting that if $p \geq 3$ then
$w \in R$ and that if $p=2$ then $x_1=3$ and at least one end vertex of the paths is in $R$. \qed

\begin{lemma}\label{SquashedJellyfish}
Suppose that there exists an $(M)^*$-packing of $K_{a,b}$ with a leave $L$ of size $l$ with $k$
non-trivial components such that exactly one vertex of $L$ has degree $4$, and every other vertex
of $L$ has degree $2$ or degree $0$. If $m_1$ and $m_2$ are even integers such that $m_1,m_2 \geq
\max(4,k+1)$ and $m_1+m_2=l$, then there is an $(M,m_1,m_2)^*$-decomposition of $K_{a,b}$.
\end{lemma}

\proof Let $\mathcal{P}$ be an $(M)^*$-packing of $K_{a,b}$ which satisfies the hypotheses of the
lemma.  By applying Lemma \ref{AttachCycles} to $\mathcal{P}$ we can obtain an $(M)^*$-packing of
$K_{a,b}$ with a leave whose only non-trivial component is a chain which has a decomposition into
an $m_1$-path and an $m_2$-path. We can then obtain an $(M,m_1,m_2)^*$-decomposition of $K_{a,b}$
by applying Lemma \ref{PathOfCycles}. \qed

\begin{lemma}\label{EquitableOneStep}
Suppose that there exists an $(M)^*$-packing of $K_{a,b}$ with a leave $L$. If $u$ and $v$ are
vertices in the same part of $K_{a,b}$ such that $\deg_L(u) > \deg_L(v)$, then there exists an
$(M)^*$-packing of $K_{a,b}$ with a leave $L'$ such that $\deg_{L'}(u)=\deg_L(u)-2$,
$\deg_{L'}(v)=\deg_L(v)+2$, and $\deg_{L'}(x)=\deg_L(x)$ for all $x \in V(L) \setminus \{u,v\}$.
Furthermore, this $L'$ also satisfies
\begin{itemize}
    \item [(i)]
if $\deg_L(v)=0$ and $u$ is not a cut vertex of $L$, then $L'$ has the same number of
non-trivial components as $L$; and
    \item [(ii)]
if $\deg_L(v)=0$, then either $L'$ has the same number of non-trivial components as $L$ or $L'$
has one more non-trivial component than $L$.
\end{itemize}
\end{lemma}

\proof Let $\mathcal{P}$ be an $(M)^*$-packing of $K_{a,b}$ which satisfies the hypotheses of the
lemma. Since $L$ is even, we have that $\deg_L(u) \geq \deg_L(v)+2$ and thus there exists a
$(u,v)$-switch $S$ whose origin is a neighbour in $L$ of $u$ and whose terminus is another
neighbour in $L$ of $u$. Let $\mathcal{P}'$ be the $(M)^*$-packing of $K_{a,b}$ obtained from
$\mathcal{P}$ by performing $S$. Using the fact that $L$ is an even graph it can be seen that
$\mathcal{P}'$ has the required properties. \qed

For an $(M)^*$-packing $\mathcal{P}$ of a graph $G$ we define
$$d(\mathcal{P})=\tfrac{1}{2}\sum_{x\in D}(\deg_{L}(x)-2),$$
where $L$ is the leave of $\mathcal{P}$ and $D$ is the set of vertices of $L$ having degree at
least $4$.

\begin{lemma}\label{FlattenOut}
Suppose that there exists an $(M)^*$-packing $\mathcal{P}_0$ of $K_{a,b}$, where $a \leq b$, with a
leave $L_0$ of size $l$, where $l \leq 2a+2$ if $a < b$ and $l \leq 2a$ if $a = b$, with $k_0$
non-trivial components such that $L_0$ has at least one vertex of degree at least $4$. Then there
exists an $(M)^*$-packing of $K_{a,b}$ with a leave $L'$ such that exactly one vertex of $L'$ has
degree $4$, every other vertex of $L'$ has degree $2$ or degree $0$, and $L'$ has at most
$\min(k_0+d(\mathcal{P}_0)-1,\lfloor\frac{l}{4}\rfloor-1)$ non-trivial components.
\end{lemma}

\proof Let $d=d(\mathcal{P}_0)$. Note that, informally, $d-1$ is the minimum number of applications
of Lemma \ref{EquitableOneStep} that would be required to transform $\mathcal{P}_0$ into an
$(M)^*$-packing of $K_{a,b}$ in whose leave one vertex has degree $4$ and every other vertex has
degree $2$ or degree $0$.

Create a sequence
$$\mathcal{P}_0,\mathcal{P}_1,\ldots,\mathcal{P}_{d-1}$$
of $(M)^*$-packings of $K_{a,b}$ inductively by letting $\mathcal{P}_{i+1}$ be the $(M)^*$-packing
of $K_{a,b}$ obtained by applying Lemma \ref{EquitableOneStep} to $\mathcal{P}_i$, choosing $u$ and
$v$ to be vertices in the same part of $K_{a,b}$ such that $u$ is of degree at least $4$ in the
leave of $\mathcal{P}_i$ and $v$ is an isolated vertex in the leave of $\mathcal{P}_i$. (To see
that such vertices exist, consider two cases depending on whether all the vertices with degree at
least $4$ in the leave of $\mathcal{P}_i$ are in the same part of $K_{a,b}$, noting that $l \leq
2a+2$ if $a < b$, that $l \leq 2a$ if $a = b$, and that for each $j \in \{0,1,\ldots,d-2\}$ the
leave of $\mathcal{P}_j$ contains at least two vertices of degree at least $4$ or at least one
vertex of degree at least $6$.) For each $i \in \{1,2,\ldots,d-1\}$, let $L_i$ be the leave of
$\mathcal{P}_i$ and let $k_i$ be the number of non-trivial components in $L_i$. Note that exactly
one vertex of $L_{d-1}$ has degree $4$ and every other vertex of $L_{d-1}$ has degree $2$ or degree
$0$.

It only remains to show that $k_{d-1} \leq \min(k_0+d-1,\lfloor\frac{l}{4}\rfloor-1)$. Since
exactly one vertex of $L_{d-1}$ has degree $4$ and every other vertex of $L_{d-1}$ has degree $2$
or degree $0$, it follows that $L_{d-1}$ has a decomposition into $k_{d-1}+1$ cycles. Since each of
these cycles has length at least $4$, we have that $l \geq 4(k_{d-1}+1)$ and so $k_{d-1} \leq
\lfloor\frac{l}{4}\rfloor-1$. Also, by Property (ii) of Lemma \ref{EquitableOneStep}, $k_{i+1} \leq
k_i + 1$ for each $i\in\{0,1,\ldots,d-2\}$ and hence $k_{d-1} \leq k_0+d-1$. \qed

\begin{lemma}\label{JoinCycles}
Suppose there exists an $(M,h,m,m')^*$-decomposition of $K_{a,b}$, where $a \leq b$. If $m+m' \leq
3h$, $m+m'+h \leq 2a+2$ if $a < b$, and $m+m'+h \leq 2a$ if $a = b$, then there exists an
$(M,h,m+m')^*$-decomposition of $K_{a,b}$.
\end{lemma}

\proof Let $\mathcal{D}$ be an $(M,h,m,m')^*$-decomposition of $K_{a,b}$ which satisfies the
hypotheses of the lemma. Let $\mathcal{P}$ be an $(M)^*$-packing of $K_{a,b}$ obtained by omitting
an $m$-cycle, an $m'$-cycle and an $h$-cycle from $\mathcal{D}$, let $L$ be the leave of
$\mathcal{P}$, and let $k$ be the number of non-trivial components of $L$. The proof now splits
into cases depending on the properties of the graph $L$.

{\bf Case 1.} Suppose that $k=3$. Then the non-trivial components of $L$ are an $m$-cycle, an
$m'$-cycle and an $h$-cycle. Let $u$ and $v$ be vertices in the same part of $K_{a,b}$ such that
$u$ and $v$ are in two distinct cycles of $L$ which have lengths $m$ and $m'$ respectively, and let
$y$ be a neighbour in $L$ of $u$. Let $\mathcal{P}^{\dag}$ be the $(M)^*$-packing of $K_{a,b}$
obtained from $\mathcal{P}$ by performing the $(u,v)$-switch $S$ with origin $y$ and let $L^{\dag}$
be the leave of $\mathcal{P}^{\dag}$. Then $L^{\dag}$ has exactly two non-trivial components, one
an $h$-cycle and the other either an $(m+m')$-cycle or an $(m,m')$-chain. In the former case we can
add the two cycles to $\mathcal{P}^{\dag}$ to complete the proof. In the latter case, since $m,m',h
\geq 4$, we can apply Lemma \ref{SquashedJellyfish} (setting $m_1=h$ and $m_2=m+m'$) to complete
the proof.

{\bf Case 2.} Suppose that $k=2$, and that the non-trivial components of $L$ are a cycle and a
$2$-chain. Then, since $m,m',h \geq 4$, we can apply Lemma \ref{SquashedJellyfish} (setting $m_1=h$
and $m_2=m+m'$) to complete the proof.

{\bf Case 3.} Suppose that we are in neither Case 1 nor Case 2. Then $L$ contains at least two
vertices of degree at least $4$ or at least one vertex of degree at least $6$. Let $d =
d(\mathcal{P})$.

We claim that it suffices to find an $(M)^*$-packing $\mathcal{P}^{\ddag}$ of $K_{a,b}$ such that
exactly one vertex in the leave of $\mathcal{P}^{\ddag}$ has degree $4$, every other vertex in the
leave of $\mathcal{P}^{\ddag}$ has degree $2$ or degree $0$, and the number $k^{\ddag}$ of
non-trivial components in the leave of $\mathcal{P}^{\ddag}$ obeys $k^{\ddag} \leq
\min(m+m'-1,\lfloor\frac{m+m'+h}{4}\rfloor-1)$. If we have such a packing then we can apply Lemma
\ref{SquashedJellyfish} (setting $m_1=h$ and $m_2=m+m'$) to it to complete the proof, provided that
$m+m'\geq k^{\ddag}+1$ and $h \geq k^{\ddag}+1$ both hold (note that $m,m',h \geq 4$). Since
$k^{\ddag} \leq m+m'-1$ we have that $m+m' \geq k^{\ddag}+1$. Since $k^{\ddag} \leq
\lfloor\frac{m+m'+h}{4}\rfloor-1$, and since $m+m' \leq 3h$ from the hypotheses of the lemma, we
have that $k^{\ddag} \leq h-1$ and thus $h \geq k^{\ddag}+1$.

So it only remains to find such an $(M)^*$-packing $\mathcal{P}^{\ddag}$. Because $L$ contains an
$h$-cycle, it has at least $h$ vertices of positive degree. Our proof now splits into cases
depending on whether $L$ has exactly $h$ vertices of positive degree or not.

{\bf Case 3a.} Suppose that $L$ has at least $h+1$ vertices of positive degree. Let
$\mathcal{P}^{\ddag}$ be the $(M)^*$-packing of $K_{a,b}$ obtained by applying Lemma
\ref{FlattenOut} to $\mathcal{P}$, let $L^{\ddag}$ be the leave of $\mathcal{P}^{\ddag}$, and let
$k^{\ddag}$ be the number of non-trivial components in $L^{\ddag}$. Note that exactly one vertex of
$L^{\ddag}$ has degree $4$, every other vertex of $L^{\ddag}$ has degree $2$ or degree $0$, and
$k^{\ddag} \leq \min(k+d-1,\lfloor\frac{m+m'+h}{4}\rfloor-1)$. Since we are not in Case 1, $k \in
\{1,2\}$. If $k=1$ then, since $L$ has at least $h+1$ vertices of positive degree and since $L$ is
an even graph with total degree $2(m+m'+h)$, it can be seen that $d \leq m+m'-1$ and hence that
$k^{\ddag} \leq \min(m+m'-1,\lfloor\frac{m+m'+h}{4}\rfloor-1)$. If $k=2$ then $L$ must have at
least $h+4$ vertices of positive degree (note that $m,m' \geq 4$), and since $L$ is an even graph
with total degree $2(m+m'+h)$, it can be seen that $d \leq m+m'-4$ and hence that $k^{\ddag} \leq
\min(m+m'-3,\lfloor\frac{m+m'+h}{4}\rfloor-1)$.

{\bf Case 3b.} Suppose that $L$ has exactly $h$ vertices of positive degree. Then, since $L$
contains an $h$-cycle, it follows that $k=1$ and that $L$ has no cut vertex. Let
$\mathcal{P}^{\dag}$ be the result of applying Lemma \ref{EquitableOneStep} to $\mathcal{P}$,
choosing $u$ and $v$ to be vertices in the same part of $K_{a,b}$ such that $u$ is of degree at
least $4$ in $L$ and $v$ is an isolated vertex in $L$. (To see that such vertices exist, consider
two cases depending on whether all the vertices with degree at least $4$ in $L$ are in the same
part of $K_{a,b}$, noting that $l \leq 2a+2$ if $a < b$, that $l \leq 2a$ if $a = b$, and that $L$
contains at least two vertices of degree at least $4$ or at least one vertex of degree at least
$6$.) Note that $d(\mathcal{P}^{\dag})=d-1$ and that, by Property (i) of Lemma
\ref{EquitableOneStep}, the leave of $\mathcal{P}^{\dag}$ has exactly one non-trivial component.
Now let $\mathcal{P}^{\ddag}$ be the $(M)^*$-packing of $K_{a,b}$ obtained by applying Lemma
\ref{FlattenOut} to $\mathcal{P^{\dag}}$, let $L^{\ddag}$ be the leave of $\mathcal{P}^{\ddag}$,
and let $k^{\ddag}$ be the number of non-trivial components in $L^{\ddag}$. Note that exactly one
vertex of $L^{\ddag}$ has degree $4$, every other vertex of $L^{\ddag}$ has degree $2$ or degree
$0$, and $k^{\ddag} \leq \min(d-1,\lfloor\frac{m+m'+h}{4}\rfloor-1)$. Since $L$ has exactly $h$
vertices of positive degree, and since $L$ is an even graph with total degree $2(m+m'+h)$, it can
be seen that $d = m+m'$ and hence that $k^{\ddag} \leq
\min(m+m'-1,\lfloor\frac{m+m'+h}{4}\rfloor-1)$. \qed

\section{Decompositions into short cycles} \label{BipartiteDecompSection}

In this section we will prove Theorem \ref{MainTheorem}. Our main tool for this task will be Lemma
\ref{JoinCycles}. We will also require the following special case of the main result of
\cite{ChFuHu}.

\begin{theorem}[\cite{ChFuHu}]\label{46Decomp}
Let $a$ and $b$ be integers such that $a,b \geq 4$ and either both $a$ and $b$ are even or $a=b$.
Let $m_1,m_2,\ldots,m_t$ be integers such that $m_1 + m_2 + \cdots + m_t = ab$ if $a$ and $b$ are
even, $m_1 + m_2 + \cdots + m_t = a^2-a$ if $a=b$ and $a$ is odd, and $m_i \in \{4,6\}$ for all $i
\in \{1,2,\ldots,t\}$. Then there exists an $(m_1,m_2,\ldots,m_t)^*$-decomposition of $K_{a,b}$.
\end{theorem}

\noindent{\bf Proof of Theorem \ref{MainTheorem}} \quad We may assume without loss of generality
that $a \leq b$. We claim that there exists an $(M)^*$-decomposition of $K_{a,b}$ for any list of
even integers $M=m_1,m_2,\ldots,m_t$ which satisfies $m_1+m_2+\cdots+m_t=a^2-a$ if $a$ is odd (and
hence $a=b$), $m_1+m_2+\cdots+m_t=ab$ if $a$ and $b$ are even, and $4\leq m_1\leq m_2 \leq \cdots
\leq m_t \leq \min(a,b,3m_{t-1})$. Proving this claim will prove the theorem.

If $a \leq 3$ then the claim is vacuous, so we may assume that $a \geq 4$. Suppose for a
contradiction that the claim is false. Let $Z=z_1,z_2,\ldots,z_t$ be a non-decreasing list which
satisfies the conditions of the claim but for which there does not exist a $(Z)^*$-decomposition of
$K_{a,b}$. Further suppose that of all such lists $Z$ is one with a maximum number of entries.

Since there does not exist a $(Z)^*$-decomposition of $K_{a,b}$, it follows from Theorem
\ref{46Decomp} that $z_t \geq 8$. Let $Z^{\dag}$ be the list $z_1, z_2, \ldots, z_{t-1}, 4, z_t-4$,
reordered so as to be non-decreasing. Since $Z$ satisfies the conditions of the claim, it is
routine to check that $Z^{\dag}$ also does. Thus, since $Z^{\dag}$ has more entries than $Z$, there
exists a $(Z^{\dag})^*$-decomposition $\mathcal{D}^{\dag}$ of $K_{a,b}$ by our definition of $Z$.
Now, provided that $z_t \leq 3z_{t-1}$ and $z_{t-1}+z_t \leq 2a$, we can apply Lemma
\ref{JoinCycles} to $\mathcal{D}^{\dag}$ (choosing $m=4$, $m'=z_t-4$ and $h=z_{t-1}$) to show that
a $(Z)^*$-decomposition of $K_{a,b}$ exists and hence obtain a contradiction. Since $Z$ satisfies
the conditions of the claim, we have that $z_t \leq \min(a,3z_{t-1})$ and hence that $z_t \leq
3z_{t-1}$ and $z_{t-1}+z_t \leq 2z_t \leq 2a$. This completes the proof. \qed

\section{Incomplete even-cycle systems}

In this section we prove Theorem \ref{HoleTheorem}. If $M$ is a list with every entry equal to some
integer $m$, then we call an $(M)$-packing of a graph $G$ an \emph{$m$-cycle packing of $G$} and we
call an $(M)$-decomposition of a graph $G$ an \emph{$m$-cycle decomposition of $G$}. Before we
prove Theorem \ref{HoleTheorem} we require two lemmas which will also prove useful in the next
section.

\begin{lemma}\label{TwoPathLeave}
Let $m \geq 4$ be an even integer, let $G$ be a complete bipartite graph each of whose parts has
even size at least $m+2$, let $R$ be a part of $G$, and let $l$ be an integer in
$\{4,6,8,\ldots,2m-4\}$ such that $|E(G)| \equiv l \mod m)$. If $p$ and $q$ are positive even
integers such that $p,q \geq l-m$ and $p+q=l$, then there is an $m$-cycle packing of $G$ whose
leave has a decomposition into a $p$-path and a $q$-path such that both end vertices of the paths
are in $R$.
\end{lemma}

\proof If $l \leq m+2$ then by Theorem \ref{MainTheorem} there is an $m$-cycle packing of $G$ with
a leave whose only non-trivial component is an $l$-cycle, and it is easy to see that this $l$-cycle
decomposes into the required paths. Thus we may suppose that $l \geq m+4$.

By Theorem \ref{MainTheorem} there is an $m$-cycle packing of $G$ with a leave whose only
non-trivial component is an $(l-m)$-cycle, and it is easy to see that there is an $m$-cycle in the
packing which shares at least one vertex with the $(l-m)$-cycle in the leave. Removing this
$m$-cycle from the packing results in an $m$-cycle packing of $G$ with a leave $L$ of size $l$ such
that $L$ has exactly one non-trivial component, $L$ has at least one and at most $l-m$ vertices of
degree 4, and every other vertex of $L$ has degree $2$ or degree $0$. Furthermore, if $L$ has
exactly $l-m$ vertices of degree $4$ then $L$ has no cut vertex, since in this case the non-trivial
component of $L$ has exactly $m$ vertices and contains an $m$-cycle. So it can be seen that, by
applying Lemma \ref{FlattenOut} to this packing, we can obtain an $m$-cycle packing of $G$ with a
leave $L'$ of size $l$ such that exactly one vertex of $L'$ has degree $4$, every other vertex of
$L'$ has degree $2$ or degree $0$, and $L'$ has at most $l-m-1$ non-trivial components. Then by
applying Lemma \ref{AttachCycles} we can obtain an $m$-cycle packing of $G$ whose leave has a
decomposition into a $p$-path and a $q$-path such that the end vertices of the paths are in $R$, as
required (note that $p,q \geq l-m \geq 4$). \qed

\begin{lemma}\label{FindSets}
Let $A$ be a set, let $S$ and $T$ be subsets of $A$ and let $s'$ and $t'$ be positive integers.
Then there exist subsets $S'$ and $T'$ of $A$ such that $|S'|=s'$, $|T'|=t'$, $S \cap S' =
\emptyset$, $T \cap T' = \emptyset$ and $S' \cap T' = \emptyset$ if and only if
\begin{itemize}
    \item[(i)]
$|S \cap T|+s'+t' \leq |A|$;
    \item[(ii)]
$|S|+s' \leq |A|$;
    \item[(iii)]
$|T|+t' \leq |A|$.
\end{itemize}
\end{lemma}

\proof To see that Conditions (i), (ii) and (iii) are necessary for the existence of $S'$ and $T'$,
note that $S \cap T$, $S'$ and $T'$ are pairwise disjoint subsets of $A$, that $S$ and $S'$ are
disjoint subsets of $A$, and that $T$ and $T'$ are disjoint subsets of $A$. We will show that these
conditions are also sufficient for the existence of the sets.

If $s' \leq |T \setminus S|$, then there exists a set $S'$ such that $|S'|=s'$ and $S' \subseteq  T
\setminus S$. Since $|T|+t' \leq |A|$, there exists a set $T'$ such that $|T'|=t'$ and $T'
\subseteq A \setminus T$. Since $S' \subseteq T \setminus S$ it follows that $S' \cap T' =
\emptyset$ and thus $S'$ and $T'$ have the required properties.

If $s' > |T \setminus S|$, then, since $|S|+s' \leq |A|$, there exists a set $S'$ such that
$|S'|=s'$, $S' \subseteq A \setminus S$ and $T \setminus S \subseteq S'$. Then, since $|S \cap
T|+s'+t' \leq |A|$, there is a set $T'$ such that $|T'|=t'$ and $T' \subseteq A \setminus ((S \cap
T) \cup S')$. Since $T \setminus S \subseteq S'$ it follows that $T \cap T' = \emptyset$ and thus
$S'$ and $T'$ have the required properties. \qed

If $A$, $S$ and $T$ are sets and $s'$ and $t'$ are integers such that $S,T \subseteq A$ and
Conditions (i), (ii) and (iii) of Lemma \ref{FindSets} hold, then we will say that the triple
$(A,S,T)$ is $(s',t')$-good.

\noindent{\bf Proof of Theorem \ref{HoleTheorem}.} Let $U$ and $W$ be disjoint sets of sizes $u$
and $v-u$ respectively. Let $\infty$ be an element of $U$. Let $G_1$ be the complete graph on
vertex set $W \cup \{\infty\}$ and let $G_2$ be the complete bipartite graph with parts $U\setminus
\{\infty\}$ and $W$. Note that $G_1 \cup G_2$ is isomorphic to $K_v - K_u$ and hence that
$|E(G_1)|+|E(G_2)| = \binom{v}{2}-\binom{u}{2}$.

Suppose first that $v-u=4$. Then we have that $m=4$ (since $4 \leq m \leq v-u$) and hence, from the
hypotheses of the theorem, that $\binom{v}{2}-\binom{u}{2} \equiv 0 \mod 4)$. Also, we have that
$|E(G_1)|=\binom{v-u+1}{2}=10$, that $|E(G_2)| \equiv 0 \mod 4)$ (since $|E(G_2)|=(u-1)(v-u)$), and
hence that $\binom{v}{2}-\binom{u}{2} \equiv 2 \mod 4)$ (since $|E(G_1)|+|E(G_2)| =
\binom{v}{2}-\binom{u}{2}$). This is a contradiction and hence we may assume that $v-u \geq 6$.

Let $e$ be the unique element of $\{0,4,6,8,\ldots,m-2,m+2\}$ such that $|E(G_1)| \equiv e \mod
m)$. It follows from Lemma 2.5 of \cite{Ho} (noting that $v-u+1 \geq 7$) that there is an $m$-cycle
packing $\mathcal{P}_1$ of $G_1$ whose leave has an $e$-cycle as its only non-trivial component if
$e \neq 0$ and whose leave is empty if $e=0$. To complete the proof we will find an $m$-cycle
packing of $G_2$ with a leave $L_2$ such that the union of $L_2$ with the leave of $\mathcal{P}_1$
(after a relabelling of the vertices in $\mathcal{P}_1$) has a decomposition into $m$-cycles. Since
$|E(G_1)|+|E(G_2)| = \binom{v}{2}-\binom{u}{2}$ and $\binom{v}{2}-\binom{u}{2} \equiv 0 \mod m)$,
we have that $|E(G_2)| \equiv -e \mod m)$. The proof now splits into cases depending on the value
of $e$.

{\bf Case 1.} Suppose that $e=0$. Then $|E(G_2)| \equiv 0 \mod m)$, by Theorem \ref{MainTheorem}
there is an $m$-cycle decomposition $\mathcal{P}_2$ of $G_2$, and $\mathcal{P}_1 \cup
\mathcal{P}_2$ is the required decomposition.

{\bf Case 2.} Suppose that $4 \leq e \leq m-2$. Since $e \neq 0$ it must be that $u \geq m+3$ and
$v-u \geq m+2$ (for otherwise either $u-1=m$ or $v-u =m$ and in either case we have that $|E(G_2)|
\equiv 0 \mod m)$ and $e=0$). By Lemma \ref{TwoPathLeave} there is an $m$-cycle packing
$\mathcal{P}_2$ of $G_2$ whose leave has a decomposition into an $(m-2)$-path $P_2$ and an
$(m-e+2)$-path $Q_2$ such that the end vertices of $P_2$ and $Q_2$ are in $W$ (note that $|E(G_2)|
\equiv 2m-e \mod m)$, that $m-2 \geq m-e$ and that $m-e+2 \geq m-e$). Let $y$ and $z$ be the
end vertices of these paths.

We claim that we can relabel the vertices in $\mathcal{P}_1$ in such a way that its leave has a
decomposition into a $2$-path $P_1$ from $y$ to $z$ and an $(e-2)$-path $Q_1$ from $y$ to $z$ such
that $V(P_1) \cap V(P_2) = \{y,z\}$ and $V(Q_1) \cap V(Q_2) = \{y,z\}$. To see that we can do this,
note that $(W,V(P_2) \cap W,V(Q_2) \cap W)$ is $(1,e-3)$-good (since $|W| \geq m+2$, $|V(P_2) \cap
W| = \frac{m}{2}$, $|V(Q_2) \cap W| = \frac{m-e+4}{2}$, and hence $|(V(P_2) \cap W) \cap (V(Q_2)
\cap W)| \leq \frac{m-e+4}{2}$) and use Lemma \ref{FindSets}. Then $\mathcal{P}_1 \cup
\mathcal{P}_2 \cup \{P_1 \cup P_2, Q_1 \cup Q_2\}$ is the required decomposition.

{\bf Case 3.} Suppose that $e = m+2$. Again, since $e \neq 0$ it must be that $u \geq m+3$ and $v-u
\geq m+2$. By Lemma \ref{TwoPathLeave} there is an $m$-cycle packing $\mathcal{P}_2$ of $G_2$ whose
leave has a decomposition into an $(m-4)$-path $P_2$ and a $2$-path $Q_2$ such that the
end vertices of the paths are in $W$ (note that $|E(G_2)| \equiv m-2 \mod m)$, that $m-4 \geq -2$
and that $2 \geq -2$). Let $y$ and $z$ be the end vertices of these paths.

We claim that we can relabel the vertices in $\mathcal{P}_1$ in such a way that its leave has a
decomposition into a $4$-path $P_1$ from $y$ to $z$ and an $(m-2)$-path $Q_1$ from $y$ to $z$ such
that $V(P_1) \cap V(P_2) = \{y,z\}$ and $V(Q_1) \cap V(Q_2) = \{y,z\}$. To see that we can do this,
note that $(W,V(P_2) \cap W,V(Q_2) \cap W)$ is $(3,m-3)$-good (since $|W| \geq m+2$, $|V(P_2) \cap
W| = \frac{m-2}{2}$, $|V(Q_2) \cap W| = 2$, and hence $|(V(P_2) \cap W) \cap (V(Q_2) \cap W)| \leq
2$) and use Lemma \ref{FindSets}. Then $\mathcal{P}_1 \cup \mathcal{P}_2 \cup \{P_1 \cup P_2, Q_1
\cup Q_2\}$ is the required decomposition. \qed

\section{Decompositions of Complete Multipartite Graphs}

In this section we will use Theorem \ref{MainTheorem} to address the problem of decomposing a
complete multipartite graph into even cycles. In particular we will prove Theorem
\ref{MultipartiteTheorem}.

\begin{lemma}\label{CombinePaths}
Let $H_2$ be a complete bipartite graph with parts $A$ and $B$ of even size, and let $H_1$ be an
even graph with vertex set $A$. If there exists an $(M_1)$-packing of $H_1$ whose leave has a
decomposition into a $p$-path $P$ and a $q$-path $Q$, and there exists an $(M_2)$-packing of $H_2$
with a leave whose only non-trivial component is an $l$-cycle, then the following hold.
\begin{itemize}
    \item[(a)]
For any even integers $p'$ and $q'$ such that $p',q' \geq 2$ and $p'+q'=l$, if $(A,V(P),V(Q))$
is $(\frac{p'-2}{2},\frac{q'-2}{2})$-good, then there exists an
$(M_1,M_2,p+p',q+q')$-decomposition of $H_1 \cup H_2$.
    \item[(b)]
For any even integers $p'$ and $q'$ such that $p',q' \geq 2$ and $p'+q'=l$, if $(A,V(P),V(Q))$
is $(\frac{p'-2}{2},\frac{q'-2}{2})$-good, then there exists an $(M_1,M_2)$-packing of $H_1
\cup H_2$ whose leave has a decomposition into a $(p+p')$-path and a $(q+q')$-path.
    \item[(c)]
For any even integers $p'$, $q'$ and $p''$ such that $p' \geq 2$, $q' \geq 0$, $p'' \geq 2$ and
$p'+q'+p''=l$, if $(A,V(P),V(Q))$ is $(\frac{p'-2}{2},\max(\frac{q'-2}{2},0))$-good, then there
exists an $(M_1,M_2,p+p')$-packing of $H_1 \cup H_2$ whose leave has a decomposition into a
$p''$-path and a $(q+q')$-path.
\end{itemize}
\end{lemma}

\proof Let $\mathcal{P}_1$ be an $(M_1)$-packing of $H_1$ whose leave has a decomposition into a
$p$-path $P$ and a $q$-path $Q$, and let $\mathcal{P}_2$ be an $(M_2)$-packing of $H_2$ with a
leave whose only non-trivial component is an $l$-cycle. Since $H_1$ is an even graph, the leave of
the $\mathcal{P}_1$ is an even graph and there are two vertices $a$ and $a^{\dag}$ in $A$ such that
$a$ and $a^{\dag}$ are the end vertices of $P$ and $Q$. We will show separately that each of (a),
(b) and (c) hold.

{\bf (a).} Using Lemma \ref{FindSets}, since $(A,V(P),V(Q))$ is
$(\frac{p'-2}{2},\frac{q'-2}{2})$-good, we can relabel vertices in $\mathcal{P}_2$ in such a way
that the leave of $\mathcal{P}_2$ has a decomposition into a $p'$-path $P'$ from $a$ to $a^{\dag}$
and a $q'$-path $Q'$ from $a$ to $a^{\dag}$ such that $V(P') \cap V(P) = V(Q') \cap V(Q) =
\{a,a^{\dag}\}$. Then $\mathcal{P}_1 \cup \mathcal{P}_2 \cup \{P \cup P',Q \cup Q'\}$ is an
$(M_1,M_2,p+p',q+q')$-decomposition of $H_1 \cup H_2$.

{\bf (b).} Let $b$ and $b^{\dag}$ be distinct vertices in $B$. Using Lemma \ref{FindSets}, since
$(A,V(P),V(Q))$ is $(\frac{p'-2}{2},\frac{q'-2}{2})$-good, we can relabel vertices in
$\mathcal{P}_2$ in such a way that the leave of $\mathcal{P}_2$ has a decomposition into the
$1$-paths $[a,b]$ and $[a,b^{\dag}]$, a $(p'-1)$-path $P'$ from $a^{\dag}$ to $b^{\dag}$ and a
$(q'-1)$-path $Q'$ from $a^{\dag}$ to $b$ such that $V(P') \cap V(P) = V(Q') \cap V(Q) =
\{a^{\dag}\}$. Then $\mathcal{P}_1 \cup \mathcal{P}_2$ is an $(M_1,M_2)$-packing of $H_1 \cup H_2$,
and $P \cup P' \cup [a,b]$ and $Q \cup Q' \cup [a,b^{\dag}]$ are the required paths in its leave.

{\bf (c).} If $q'=0$, then using Lemma \ref{FindSets}, since $(A,V(P),V(Q))$ is
$(\frac{p'-2}{2},0)$-good, we can relabel vertices in $\mathcal{P}_2$ in such a way that the leave
of $\mathcal{P}_2$ has a decomposition into a $p'$-path $P'$ from $a$ to $a^{\dag}$ and a
$p''$-path $P''$ from $a$ to $a^{\dag}$ such that $V(P') \cap V(P) = \{a,a^{\dag}\}$. Then
$\mathcal{P}_1 \cup \mathcal{P}_2 \cup \{P \cup P'\}$ is an $(M_1,M_2,p+p')$-packing of $H_1 \cup
H_2$ and $P''$ and $Q$ are the required paths in its leave.

If $q' \geq 2$, then let $b$ and $b^{\dag}$ be distinct vertices in $B$. Using Lemma
\ref{FindSets}, since $(A,V(P),V(Q))$ is $(\frac{p'-2}{2},\frac{q'-2}{2})$-good, we can relabel
vertices in $\mathcal{P}_2$ in such a way that the leave of $\mathcal{P}_2$ has a decomposition
into the $1$-path $[a,b^{\dag}]$, a $p'$-path $P'$ from $a$ to $a^{\dag}$, a $(q'-1)$-path $Q'$
from $a^{\dag}$ to $b$, and a $p''$-path $P''$ from $b$ to $b^{\dag}$ such that $V(P') \cap V(P) =
\{a,a^{\dag}\}$ and $V(Q') \cap V(Q) = \{a^{\dag}\}$. Then $\mathcal{P}_1 \cup \mathcal{P}_2 \cup
\{P \cup P'\}$ is an $(M_1,M_2,p+p')$-packing of $H_1 \cup H_2$, and $P''$ and $Q \cup Q' \cup
[a,b^{\dag}]$ are the required paths in its leave. \qed

By combining Theorem \ref{MainTheorem} and Lemma \ref{CombinePaths}, many decompositions of
complete multipartite graphs into even length cycles can be created. Here, we will restrict our
focus to decompositions of complete multipartite graphs into uniform even length cycles and content
ourselves with proving Theorem \ref{MultipartiteTheorem}. The proof of this result is indicative of
how Theorem \ref{MainTheorem} and Lemma \ref{CombinePaths} can be combined effectively.

\noindent{\bf Proof of Theorem \ref{MultipartiteTheorem}.} If $G$ has only one part then it has no
edges and the result is trivial, so we may assume that $G$ has at least two parts. In \cite{CaBi}
necessary and sufficient conditions are found for the existence of a decomposition of a complete
multipartite graph into $m$-cycles in the cases $m=4$ and $m=6$. It follows from these conditions
that the theorem is true for $m \in \{4,6\}$, so we may assume that $m \geq 8$. Let the parts of
$G$ be $V_1,V_2,\ldots,V_t$. For each $i \in \{1,2,\ldots,t-1\}$ let $G_i$ be the complete
bipartite graph with parts $V_1\cup V_2\cup \cdots\cup V_i$ and $V_{i+1}$, and note that $G$ is the
edge disjoint union of the graphs $G_1, G_2, \ldots, G_{t-1}$. Also, for each $i \in
\{1,2,\ldots,t-1\}$, let $l_i$ be the element of $\{0,4,6,8,\ldots,m-2,m+2\}$ such that $|E(G_i)|
\equiv l_i \mod m)$.

Roughly speaking, we will find the required $m$-cycle decomposition of $G$ inductively, by
beginning with an $m$-cycle packing of $G_1$ (obtained using Theorem \ref{MainTheorem}) and, at the
$i$th stage, combining the $m$-cycle packing of $G_1 \cup G_2 \cup \cdots \cup G_i$ which we have
already created with an $m$-cycle packing of $G_{i+1}$ whose leave has size $l_{i+1}$ (this packing
of $G_{i+1}$ is also obtained using Theorem \ref{MainTheorem}). The packings are combined by either
simply taking their union or by applying Lemma \ref{CombinePaths}.

Let $\sigma_i = l_{i+1}+l_{i+2}+\cdots+l_{t-1}$ for each $i \in \{1,2,\ldots,t-2\}$, and let
$\sigma_{t-1} = 0$. Intuitively, $\sigma_i$ is the sum of the sizes of the leaves of the packings
which we are yet to add, when we have just constructed our packing of $G_1 \cup G_2 \cup \cdots
\cup G_i$.

We claim that for each $s \in \{1,2,\ldots,t-1\}$ there is an $m$-cycle packing $\mathcal{P}_s$ of
$G_1 \cup G_2 \cup \cdots \cup G_s$ with a leave $L_s$ such that
\begin{itemize}
    \item[(i)]
if $\sigma_s \equiv 0 \mod 2m)$, then $L_s$ is empty;
    \item[(ii)]
if $\sigma_s \equiv 2 \mod 2m)$, then $L_s$ has a decomposition into an $(m-4)$-path and a
$2$-path;
    \item[(iii)]
if $\sigma_s \equiv 2m-2 \mod 2m)$, then $L_s$ has a decomposition into an $(m-2)$-path and a
$4$-path;
    \item[(iv)]
if $\sigma_s \equiv x \mod 2m)$ for some $x \in \{4,6,8,\ldots,2m-4\}$, then $L_s$ has a
decomposition into a $d_s$-path and an $e_s$-path, where $d_s$ and $e_s$ are the even integers
such that $d_s+e_s=2m-x$ and $d_s \in \{e_s,e_s+2\}$.
\end{itemize}
Since $\sigma_{t-1}=0$, when $s=t-1$ this claim is the result of the theorem, and hence it suffices
to prove the claim. We shall do so by induction on $s$. Note that, since $|E(G)| \equiv 0 \mod m)$,
we have that $l_1+l_2+\cdots+l_{t-1} \equiv 0 \mod m)$, and hence that $l_1+l_2+\cdots+l_i+\sigma_i
\equiv 0 \mod m)$ for each $i \in \{1,2,\ldots,t-1\}$.

Suppose that $s=1$. Note that, since $l_1+\sigma_1 \equiv 0 \mod m)$, it follows from the
definition of $l_1$ that $|E(G_1)| \equiv -\sigma_1 \mod m)$.  If $\sigma_1 \equiv 0 \mod 2m)$,
then $|E(G_1)| \equiv 0 \mod m)$ and by Theorem \ref{MainTheorem} there is an $m$-cycle
decomposition of $G_1$. If $\sigma_1 \equiv 2 \mod 2m)$, then $|E(G_1)| \equiv m-2 \mod m)$ and by
Theorem \ref{MainTheorem} there is an $m$-cycle packing of $G_1$ with a leave whose only
non-trivial component is an $(m-2)$-cycle. This $(m-2)$-cycle can easily be decomposed into an
$(m-4)$-path and a $2$-path. If $\sigma_1 \equiv 2m-2 \mod 2m)$, then $|E(G_1)| \equiv 2 \mod m)$
and by Theorem \ref{MainTheorem} there is an $m$-cycle packing of $G_1$ with a leave whose only
non-trivial component is an $(m+2)$-cycle. This $(m+2)$-cycle can easily be decomposed into an
$(m-2)$-path and a $4$-path. Finally, if $\sigma_1 \equiv x \mod 2m)$ for some $x \in
\{4,6,8,\ldots,2m-4\}$, then $|E(G_1)| \equiv -x \mod m)$ and by Lemma \ref{TwoPathLeave} (setting
$l=2m-x)$ it can be seen that there is an $m$-cycle packing of $G_1$ whose leave has a
decomposition into a $d_1$-path and an $e_1$-path (it is easy to see that $d_1,e_1 \geq m-x$ from
their definition). So the claim is indeed true in the case $s=1$.

Now suppose that the claim is true for $s=k$ for some integer $k$ such that $1 \leq k \leq t-2$, so
there is an $m$-cycle packing $\mathcal{P}_k$ of $G_1 \cup G_2 \cup \cdots \cup G_k$ with a leave
$L_k$ satisfying (i), (ii), (iii) and (iv) for $s=k$. We will show that the claim is true for
$s=k+1$.

Note that by Theorem \ref{MainTheorem} there is an $m$-cycle packing $\mathcal{C}$ of $G_{k+1}$
whose leave has an $l_{k+1}$-cycle as its only non-trivial component if $l_{k+1} \neq 0$ and whose
leave is empty if $l_{k+1} = 0$. If $l_{k+1} \geq 8$ then, also by Theorem \ref{MainTheorem}, for
each $z \in \{4,6,8,\ldots,l_{k+1}-4\}$ there exists a packing $\mathcal{C}^z$ of $G_k$ with some
number of $m$-cycles and one $(l_{k+1}-z)$-cycle with a leave whose only non-trivial component is a
$z$-cycle. By combining $\mathcal{P}_k$ with one of these packings, either by simply taking a union
or by applying Lemma \ref{CombinePaths} (setting $H_2=G_{k+1}$, $A=V_1\cup V_2\cup \cdots\cup
V_{k+1}$, $B=V_{k+2}$ and $H_1=G_1 \cup G_2 \cup \cdots \cup G_k$), we can obtain an $m$-cycle
packing $\mathcal{P}_{k+1}$ of $G_1 \cup G_2 \cup \cdots \cup G_{k+1}$ satisfying (i), (ii), (iii)
and (iv) for $s=k+1$. The following list gives an outline of how to accomplish this. Note that,
from our definition of $\sigma_1, \sigma_2, \ldots, \sigma_t$ we have that
$\sigma_{k+1}=\sigma_k-l_{k+1}$. Also note that in each case below where we apply Lemma
\ref{CombinePaths} its hypotheses are satisfied since $l_{k+1} \leq m+2$, the two paths in the
leave of $\mathcal{P}_k$ each have length at most $m-2$, and $|V_1 \cup V_2 \cup \cdots \cup
V_{k+1}| \geq (k+1)(m+2) \geq 2(m+2)$.
\begin{itemize}
    \item
If $l_{k+1}=0$ then $\sigma_{k+1}=\sigma_{k}$ and we take the union of $\mathcal{P}_k$ and
$\mathcal{C}$.
    \item
If $\sigma_k \equiv 0 \mod 2m)$ and $l_{k+1} \in \{4,6,8,\ldots,m-2,m+2\}$, then $\sigma_{k+1}
\equiv y \mod 2m)$ for some $y \in \{4,6,8,\ldots,2m-4\}$ (note that $m \geq 8$). We take the
union of $\mathcal{P}_k$ and $\mathcal{C}$ and note that the $l_{k+1}$-cycle in the leave of
$\mathcal{C}$ has a decomposition into a $d_{k+1}$-path and an $e_{k+1}$-path.
    \item
If $\sigma_k \equiv 2 \mod 2m)$ and $l_{k+1} = 4$, then $\sigma_{k+1} \equiv 2m-2 \mod 2m)$. We
apply Lemma \ref{CombinePaths} (b) (setting $p=m-4$, $p'=2$, $q=2$ and $q'=2$) to
$\mathcal{P}_k$ and $\mathcal{C}$.
    \item
If $\sigma_k \equiv 2 \mod 2m)$ and $l_{k+1} \in \{6,8,10,\ldots,m-2,m+2\}$, then $\sigma_{k+1}
\equiv y \mod 2m)$ for some $y \in \{4,6,8,\ldots,2m-4\}$ (note that $m \geq 8$). We apply
Lemma \ref{CombinePaths} (c) (setting $p=m-4$, $p'=4$, $q=2$, $q'=d_{k+1}-2$ and $p''=e_{k+1}$)
to $\mathcal{P}_k$ and $\mathcal{C}$.
    \item
If $\sigma_k \equiv 2m-2 \mod 2m)$ and $l_{k+1} \in \{4,6,8,\ldots,m-2,m+2\}$, then
$\sigma_{k+1} \equiv y \mod 2m)$ for some $y \in \{4,6,8,\ldots,2m-4\}$ (note that $m \geq 8$).
Note that, in fact, $y \leq 2m-6$ and hence $d_{k+1} \geq 4$. We apply Lemma \ref{CombinePaths}
(c) (setting $p=m-2$, $p'=2$, $q=4$, $q'=d_{k+1}-4$ and $p''=e_{k+1}$) to $\mathcal{P}_k$ and
$\mathcal{C}$.
    \item
If $\sigma_k \equiv x \mod 2m)$ for some $x \in \{4,6,8,\ldots,2m-4\}$ and $4 \leq l_{k+1} \leq
x-4$, then $\sigma_{k+1} \equiv y \mod 2m)$ for some $y \in \{4,6,8,\ldots,2m-4\}$. Note that,
in fact, $y \leq x-4$ and hence $d_{k+1} \geq d_k+2$ and $e_{k+1} \geq e_k+2$. We apply Lemma
\ref{CombinePaths} (b) (setting $p=d_k$, $p'=d_{k+1}-d_k$, $q=e_k$ and $q'=e_{k+1}-e_k$) to
$\mathcal{P}_k$ and $\mathcal{C}$.
    \item
If $\sigma_k \equiv x \mod 2m)$ for some $x \in \{4,6,8,\ldots,2m-4\}$ and $l_{k+1} = x-2$,
then $\sigma_{k+1} \equiv 2 \mod 2m)$. Note that, since $l_{k+1} = x-2$ and $l_{k+1} \geq 4$,
we have that $x \geq 6$ and hence $e_{k} \leq m-4$. We apply Lemma \ref{CombinePaths} (c)
(setting $p=d_k$, $p'=m-d_k$, $q=e_k$, $q'=m-4-e_k$ and $p''=2$) to $\mathcal{P}_k$ and
$\mathcal{C}$.
    \item
If $\sigma_k \equiv x \mod 2m)$ for some $x \in \{4,6,8,\ldots,2m-4\}$ and $l_{k+1} = x$, then
$\sigma_{k+1} \equiv 0 \mod 2m)$. We apply Lemma \ref{CombinePaths} (a) (setting $p=d_k$,
$p'=m-d_k$, $q=e_k$ and $q'=m-e_k$) to $\mathcal{P}_k$ and $\mathcal{C}$.
    \item
If $\sigma_k \equiv x \mod 2m)$ for some $x \in \{4,6,8,\ldots,2m-4\}$ and $l_{k+1} = x+2$,
then $\sigma_{k+1} \equiv 2m-2 \mod 2m)$. We apply Lemma \ref{CombinePaths} (c) (setting
$p=d_k$, $p'=m-d_k$, $q=e_k$, $q'=m-2-e_k$ and $p''=4$) to $\mathcal{P}_k$ and $\mathcal{C}$.
    \item
If $\sigma_k \equiv x \mod 2m)$ for some $x \in \{4,6,8,\ldots,2m-4\}$ and $x+4 \leq l_{k+1}
\leq m+2$, then $\sigma_{k+1} \equiv y \mod 2m)$ for some $y \in \{4,6,8,\ldots,2m-4\}$ (since
$m \geq 8$). Note that, in fact, $y = 2m-(l_{k+1}-x)$ and that $l_{k+1} \geq 8$. We apply Lemma
\ref{CombinePaths} (a) (setting $p=d_k$, $p'=m-d_k$, $q=e_k$ and $q'=m-e_k$) to $\mathcal{P}_k$
and $\mathcal{C}^x$, remove the $(l_{k+1}-x)$-cycle from this packing, and note that this cycle
has a decomposition into a $d_{k+1}$-path and an $e_{k+1}$-path. \qed
\end{itemize}

\vspace{0.3cm} \noindent{\bf Acknowledgement}

This research was supported in part by the Atlantic Association for Research in the Mathematical
Sciences.


\begin{thebibliography}{99}
\def\baselinestretch{1}\small\normalsize
\begin{small}

\bibitem{Al} B. Alspach, Research problems, Problem 3, {\it
    Discrete Math.} {\bf 36} (1981), 333.

\bibitem{AlGa} B. Alspach and H. Gavlas, Cycle decompositions
    of $K_n$ and $K_n - I$, {\it J.
    Combin. Theory Ser. B} {\bf 81} (2001), 77--99.

\bibitem{ArEtAl} D. Archdeacon, M. Debowsky, J. Dinitz and H. Gavlas, Cycle systems in the complete
    bipartite graph minus a one-factor, {\it Discrete Math.} {\bf 284} (2004), 37--43.

\bibitem{Bi} E.J. Billington, Multipartite graph decomposition: cycles and closed
    trails, {\it Matematiche (Catania)} {\bf 59} (2004), 53--72.

\bibitem{Br} D. Bryant, Cycle decompositions of complete
    graphs, in {\it Surveys in Combinatorics 2007},
    A. Hilton and J. Talbot (Editors),
    London Mathematical Society Lecture Note Series {\bf 346},
    Proceedings of the 21st British Combinatorial Conference, Cambridge University Press, 2007, pp.
    67--97.


\bibitem{BrHoLong} D. Bryant and D. Horsley, Decompositions of complete graphs into long cycles,
    {\it Bull. Lond. Math. Soc.} {\bf 41} (2009), 927--934.

\bibitem{BrHoBigN} D. Bryant and D. Horsley, An asymptotic solution to the cycle decomposition
    problem for complete graphs, {\it J. Combin. Theory Ser. A} {\bf 117} (2010), 1258-1284.

\bibitem{BrHoMa} D. Bryant, D. Horsley and B. Maenhaut,
    Decompositions into $2$-regular subgraphs
    and equitable partial cycle decompositions, {\it J. Combin. Theory Ser. B} {\bf 93} (2005),
    67--72.

\bibitem{BrRo} D. Bryant and C.A. Rodger, Cycle decompositions, in {\it The CRC Handbook of
    Combinatorial Designs, {\rm 2}nd edition} (Eds. C. J. Colbourn, J. H. Dinitz), CRC Press, Boca
    Raton (2007), 373--382.

\bibitem{BrRoSp} D. Bryant, C.A. Rodger and E. Spicer, Embeddings of $m$-cycle systems and
    incomplete $m$-cycle systems: $m\leq 14$, {\it Discrete Math.} {\bf 171} (1997), 55--75.

\bibitem{CaBi} N.J. Cavenagh and E.J. Billington, Decomposition of complete multipartite graphs
    into cycles of even length, {\it Graphs Combin.} {\bf 16} (2000), 49--65.

\bibitem{ChFu} C-C. Chou, C-M. Fu, Decomposition of $K_{m,n}$ into $4$-cycles and
    $2t$-cycles, {\it J. Comb. Optim.} {\bf 14} (2007), 205--218.

\bibitem{ChFuHu} C-C. Chou, C-M. Fu and W-C. Huang, Decomposition of $K_{m,n}$ into short cycles,
    {\it Discrete Math.} {\bf 197/198} (1999), 195--203.

\bibitem{HoLiRo} D.G. Hoffman, C.C. Lindner and C.A. Rodger, On the construction of odd cycle
    systems, {\it J. Graph Theory} {\bf 13} (1989), 417--426.

\bibitem{Ho} D. Horsley, Maximum packings of the complete graph with uniform length cycles, {\it J.
    Graph Theory}, (to appear).

\bibitem{LaAu} R. Laskar and B. Auerbach, On decomposition of $r$-partite graphs into
    edge-disjoint Hamilton circuits, {\it Discrete Math.} {\bf 14} (1976), 265--268.

\bibitem{MaPuSh}J. Ma, L. Pu and H. Shen, Cycle decompositions of $K\sb{n,n}-I$, {\it SIAM
    J. Discrete Math.} {\bf 20} (2006), 603--609.

\bibitem{Sa1} M. \v{S}ajna, On decomposing $K_n-I$ into cycles of a fixed odd length, {\it Discrete
    Math.} {\bf 244} (2002), 435--444.

\bibitem{Sa2} M. \v{S}ajna, Cycle decompositions III: complete graphs and fixed length cycles, {\it
    J. Combin. Des.} {\bf 10} (2002), 27--78.

\bibitem{So} D. Sotteau, Decomposition of $K_{m,n}$ ($K^*_{m,n}$) into cycles (circuits) of length
    $2k$, {\it J. Combin. Theory Ser. B} {\bf  30} (1981), 75--81.
\end{small}
\end{thebibliography}
\end{document}